\documentclass[12pt]{article}
\usepackage{amssymb}

\author{I. G. Korepanov\\
\normalsize Southern Ural State University\\[-0.5ex]
\normalsize 76 Lenin avenue\\[-0.5ex]
\normalsize 454080 Chelyabinsk, Russia\\[-0.5ex]
\normalsize E-mail: kig@susu.ac.ru}
\date{}
\title{Euclidean 4-simplices and invariants of four-dimensional
manifolds: I.~Moves~$3\to 3$}

\def\be{\begin{equation}}
\def\ee{\end{equation}}

\def\ol{\overline}
\def\pa#1#2{{\partial#1\over\partial#2}}

\def\rank{\mathop{\rm rank}\nolimits}

\newtheorem{theorem}{Theorem}

\begin{document}
\maketitle

\begin{abstract}
We construct invariants of four-dimensional piecewise-linear manifolds,
represented as simplicial complexes, with respect to rebuildings
that transform a cluster of three 4-simplices having a common
two-dimensional face in a different cluster of the same type and having
the same boundary. Our construction is based on the usage of Euclidean
geometric values.
\end{abstract}


\section{Introduction}

The present work is planned to be the first one in a series of papers
whose aim is the construction of an invariant of four-dimensional
piecewise-linear manifolds similar to the invariants of
three-dimensional manifolds constructed in papers
\cite{3dim1} and~\cite{3dim2}. We recall that in those papers we were
considering orientable manifolds represented as simplicial complexes
(or ``pre-complexes''
where the boundary of a simplex may contain several identical components).
We assigned a Euclidean length to every edge, and a sign $+$ or $-$
to every 3-cell (tetrahedron) in such way that the following ``zero
curvature'' condition was satisfied: the algebraic sum of dihedral angles
abutting on each edge was zero modulo~$2\pi$. Here each dihedral angle
is determined by the lengths of edges of the corresponding tetrahedron
and is taken with the sign that we have assigned to that tetrahedron.
Then we considered infinitesimal variations
of edge lengths and infinitesimal ``deficit angles'' at edges depending on
those length variations. The matrix that expressed these (linear) dependencies
played the key r\^ole in the resulting explicit formula for
manifold invariants.

The paper~\cite{3dim2} shows that the invariants constructed in this way
are sensitive enough in the three-dimensional case: they distinguish
(in particular) the homotopy equivalent lens spaces $L(7,1)$
and~$L(7,2)$. On the other hand, a formula derived in
preprint~\cite{korepanov 4s} suggests the possibility of
generalizing these constructions to
the four-dimensional case and, most likely, higher dimensions.
That formula deals with the re-building of {\em one\/} cluster
of three 4-simplices, i.e.\ is {\em local\/} in a natural sense.
The present work is devoted to constructing invariants of any sequences
of such $3\to 3$ re-buildings (three 4-simplices are replaced with three
different 4-simplices with the same boundary) in a ``large'' simplicial
complex.

Now we describe the contents of the following sections.
In Section~\ref{sec local} we derive the ``local'' formula, that is, the
formula for one re-building. In comparison with paper~\cite{korepanov 4s},
the formula is derived with some refinements due to the fact that we now
take into account the orientations of 4-simplices and, consequently,
the signs of dihedral angles. In Section~\ref{sec differentials} we introduce
the basic values --- lengths, areas and deficit angles of two types --- from
which our invariants will be built in the global case. In Section~\ref{sec matrices}
we prove some useful properties of matrices that express the relations
between the differentials of those values. Here, a key r\^ole is played by
the symmetry property of one of those matrices, which is deduced with the help
of the four-dimensional version of Schl\"afli differential identity
(see~\cite{milnor}). In Section~\ref{sec invariant 3-3} we obtain our
``global'' invariant of moves $3\to 3$. To do that, we need to use the
local formula in combimation with some differential forms.
In the concluding Section~\ref{sec discussion} we discuss our results.

\section{The local formula}
\label{sec local}

Consider a Euclidean 4-simplex $ABCDE$ with all edge lengths fixed
except the length of edge~$AE$. Opposite $AE$, there is a two-dimensional
face~$BCD$ where two three-dimensional faces come together,
forming an inner dihedral angle~$\vartheta_{BCD}$.
The starting point of our constructions is the formula
\be
{dL_{AE}\over V_{ABCDE}} = 24 \, {d\vartheta_{BCD}\over S_{BCD}},
\label{basic}
\ee
where we denoted by $L_{AE}$ the {\em squared\/} length~$AE$, while $V_{ABCDE}$
and $S_{BCD}$ are the corresponding hypervolume and area
(cf.~\cite[formula~(15)]{multidim}).

Now we consider six points $A,\ldots,F$ in $\mathbb R^4$ with all
distances between them fixed except $AE$ and~$AF$. Then
(cf.~\cite[(16)]{multidim})
\be
{dL_{AE}\over V_{ABCDE}} = {dL_{AF}\over V_{ABCDF}}.
\label{basic1}
\ee
Proof: the angles at the face $BCD$ in 4-simplices $ABCDE$ and $ABCDF$
differ by the value of the angle at the same face in a ``rigid'' 4-simplex $BCDEF$ where all lengths
are fixed. Formula (\ref{basic1}) follows from the equalness of length
differentials, with (\ref{basic}) taken into account.

{\bf Important remark: }\vadjust{\nobreak}in formula (\ref{basic1}),
one has to regard $V_{ABCDE}$ and $V_{ABCDF}$ as {\em oriented volumes}
and to take the dihedral angles therein with the signs coinciding with
the signs of volumes (to be exact, we assume that they lie in the interval $(0,\pi)$
for a positive volume and in $(-\pi,0)$ for a negative one).

It is clear that we can obtain relations between $dL$ for any two edges
sharing a common vertex by permutations of letters in~(\ref{basic1}).
If, however, we want to relate, e.g., $dL_{AB}$ with $dL_{DE}$ assuming
the rest of the lengths constant, we can act as follows: consider $L_{AD}$ as a function
of $L_{AB}$ and $L_{DE}$ (the other lengths are fixed) and find
$$
{dL_{DE}\over dL_{AB}} = -{\partial L_{AD}/ \partial L_{AB} \over
\partial L_{AD}/ \partial L_{DE}}
$$
as the derivative of an implicit function $L_{AD}(L_{AB},L_{DE})={\rm const}$.
In this way we find
\be
{dL_{DE} \over V_{\hat A} V_{\hat B}} = -{dL_{AB}
\over V_{\hat D} V_{\hat E}},
\label{basic2}
\ee
where we adopt the following notations: $\hat A$, $\hat B$, etc.\ is the
sequence $ABCDEF$ with a removed letter $A$, $B$, etc.;
$V_{\ldots}$ means oriented four-dimensional volume.

Consider now three 4-simplices $ABCEF=\hat D$, $ABCFD=-\hat E$ and
$ABCDE=\hat F$ with the common two-dimensional face~$ABC$. We denote
the angles in those respective simplices at that face as
$\alpha$, $\beta$ and~$\gamma$. They form together the
full angle: $\alpha+\beta+\gamma=0 \pmod{2\pi}$. Assuming that everything
is taking place in a Euclidean space and only $AB$ and $DE$ of all lengths
can change, we find:
\be
0=d(\alpha+\beta+\gamma) = \pa{\gamma}{L_{DE}}\,dL_{DE} +
\pa{(\alpha+\beta+\gamma)}{L_{AB}}\,dL_{AB}.
\label{basic3}
\ee
Here, of course, $\partial\alpha/ \partial L_{AB}$ is calculated from
simplex $\hat D$, while $\partial\beta/ \partial L_{AB}$ --- from $\hat E$ and
$\partial\gamma/ \partial L_{AB}$ --- from~$\hat F$.

We denote the sum $\alpha+\beta+\gamma$, regarded as a function of ten
formally independent squared lengths, as $-\omega_{ABC}$ (cf.~below
formula~(\ref{defect angle})). It can be regarded this way because
each of angles $\alpha,\beta,\gamma$ is calculated separately from its own
4-simplex. Clearly, all simplices can be placed together in
$\mathbb R^4$ only if $\omega_{ABC}=0 \pmod{2\pi}$.

Taking into account that
$$
\pa{\gamma}{L_{DE}} = {1\over 24} {S_{ABC}\over V_{\hat F}}
$$
(cf.~(\ref{basic})) and comparing (\ref{basic3}) with (\ref{basic2}),
we get:
\be
\pa{\omega_{ABC}}{L_{AB}} = -{S_{ABC}\over 24}
{V_{\hat A} V_{\hat B} \over V_{\hat D} V_{\hat E} V_{\hat F}}.
\label{basic4}
\ee

Similarly, we define $\omega_{DEF}$ as minus sum of the dihedral angles
at face $DEF$ in 4-simplices $\hat A$, $-\hat B$ and $\hat C$, and
get the equality
\be
\pa{\omega_{DEF}}{L_{DE}} = -{S_{DEF}\over 24}
{V_{\hat D} V_{\hat E} \over V_{\hat A} V_{\hat B} V_{\hat C}}.
\label{basic5}
\ee
The value $\omega_{DEF}$ can, too, be considered as a function of
ten formally independent~$L$.

We will be dealing only with infinitesimal deviations of values~$L$ from
the ``flat'' case $\omega_{ABC}=\omega_{DEF}=0$. In such case, the
differentials $d\omega_{ABC}$ and $d\omega_{DEF}$ are {\em proportional},
because they are two linear forms, both depending on ten $dL$, whose
null spaces (kernels) coincide (because both conditions $d\omega_{ABC}=0$
and $d\omega_{DEF}=0$ mean that the corresponding infinitesimal
deformation of edge lengths can be realized in~$\mathbb R^4$).
In particular, if only $L_{AB}$ and $L_{DE}$ can change then
\be
\matrix{
d\omega_{ABC} = c_1 \, dL_{AB} + c_2\, dL_{DE},\cr
d\omega_{DEF} = c_3 \, dL_{AB} + c_4\, dL_{DE}, }
\label{basic6}
\ee
where the equal ratios $c_1/c_2$ and $c_3/c_4$ are found by comparing
(\ref{basic6}) with the relation (\ref{basic2}) that was obtained, of course,
under the condition $d\omega_{ABC}=d\omega_{DEF}=0$. Namely,
$$
{c_1\over c_2} = {c_3\over c_4} = {V_{\hat A} V_{\hat B} \over
V_{\hat D} V_{\hat E} }.
$$
Besides, $c_1$ and $c_4$ are right-hand sides of (\ref{basic4}) and
(\ref{basic5}), respectively.

All this together allows us to calculate the proportionality factor
between $d\omega_{ABC}$ and $d\omega_{DEF}$. We write the result in the
following form:
\be
V_{\hat D} V_{-\hat E} V_{\hat F}\, {d\omega_{ABC}\over S_{ABC}} =
V_{\hat A} V_{-\hat B} V_{\hat C}\, {d\omega_{DEF}\over S_{DEF}}.
\label{6term}
\ee

The most remarkable thing in relation (\ref{6term}) is that all the values
entering in its l.h.s.\ belong to 4-simplices $\hat D$, $\hat E$ and
$\hat F$, while those in its r.h.s.\ --- to $\hat A$, $\hat B$ and~$\hat C$.
Thus, (\ref{6term}) is an algebraic relation corresponding to the
move~$3\to 3$.

We make also the following remarks. Of course, $V_{-\hat E}$ and $V_{-\hat B}$
are nothing but $-V_{\hat E}$ and~$-V_{\hat B}$. The way we have written the
formula (\ref{6term}) is due to the fact that 4-simplices
$\hat D$, $-\hat E$ and $\hat F$ have consistent orientations, and
the same applies to $\hat A$, $-\hat B$ and~$\hat C$. On the other hand,
the common boundary of 4-simplices $\hat A$, $-\hat B$ and~$\hat C$ coincides
with that of $\hat D$, $-\hat E$ and~$\hat F$ {\em with the orientation
taken into account}. The sign of each of dihedral angles entering in
$-\omega_{ABC}$ and $-\omega_{DEF}$ coincides with the sign of the respective
quantity $V_{\hat A}$, $V_{-\hat B}$, $V_{\hat C}$, $V_{\hat D}$,
$V_{-\hat E}$, or~$V_{\hat F}$.

\section{The global case: metric values and their differentials}
\label{sec differentials}

In this Section, we introduce the principal quantities --- lengths,
areas, and two types of deficit angles ---
that will take part in the game in the global case.
Relations between variations, i.e.\ differentials, of those values
will play the main r\^ole.

Consider an orientable four-dimensional manifold~$M$. We represent it
as a simplicial complex. Choose a consistent orientation of the
4-simplices in this complex (by orientation of a simplex we understand
an ordering of its vertices taken up to even permutations).

We assign lenths to the edges of the complex and signs to its 4-simplices
in full analogy with the three-dimensional case described in the Introduction
to paper~\cite{3dim2}. Namely, we consider the {\em universal
cover\/} of the triangulation of manifold~$M$ which is itself, of course, a simplicial complex.
Its vertices are divided in classes which are inverse
images of each given vertex w.r.t.\ the covering map.
Let a homomorphism (representation) $\varphi\colon\; \pi_1(M)\to E_4$ be fixed
of the fundamental group of manifold~$M$ in the group of motions of the
four-dimensional Euclidean space. We construct a mapping~$f$ from the set of
vertices of the universal cover to $\mathbb R^4$ in the following way.
Map one arbitrarily chosen vertex in each class to an arbitrary
point. Any of the remaining vertices, say, vertex~$B$,
can be represented as $B=g(A)$ for some element
$g\in \pi_1(M)$ and some vertex~$A$ whose image has already been constructed.
Demand then that $f(B) = \varphi(g)\,f(A)$.

Having constructed the images of the vertices of our complex universal cover
in~$\mathbb R^4$, we can also put in correspondence to every cell of nonzero
dimension --- an edge, a triangle, a tetrahedron, a 4-simplex --- its image
in~$\mathbb R^4$. That will be, of course, the convex hull of its vertex images.
We emphasize that we are not afraid of possible intersections of the resulting
simplices in~$\mathbb R^4$ but we demand that all the 4-simplices in
$\mathbb R^4$ --- the images of 4-simplices from the universal cover ---
be {\em nondegenerate\/} (have nonzero volume).

Thus we get lengths associated to the edges of
the universal cover --- they are, of course, the distances between the
respective vertices in~$\mathbb R^4$. It is clear that if two given
edges are mapped into the same one by the covering map then their lengths
obtained from our construction are equal (because their images
in~$\mathbb R^4$ are taken one into the other by some element
$\varphi(g)\in E_4$).

We assume that we have fixed an orientation of~$\mathbb R^4$. Then a
4-simplex from the triangulation (recall that we have chosen a
consistent orientation for them) either preserves its orientation
under the mapping in~$\mathbb R^4$, or changes it. In the first case,
we assign to it the sign ``$+$'', in the second --- the sign~``$-$''.
Clearly, if two 4-simplices are sent to the same simplex
by the covering map then they get the same sign.

This means that we have assigned lengths to the edges and signs to the
4-simplices of the manifold~$M$'s triangulation itself and not only
of its universal cover.

As soon as the edges of each 4-simplex have acquired Euclidean lengths,
one can calculate (inner) {\em dihedral angles\/} formed by pairs of
three-dimensional faces of the 4-simplex intersecting at a two-dimensional
face. Below, we consider every such angle with the sign $+$ or $-$
coinciding with the one assigned to the 4-simplex.
That is what we mean when we speak below of {\em algebraic sums\/}
of such angles.

We call the {\em deficit angle $\omega_{ABC}$\/} at a two-dimensional
face $ABC$ of the complex the algebraic sum of all dihedral angles around
that face taken with the minus sign and considered to within the
multiples of~$2\pi$:
\be
\omega_{ABC} \stackrel{\rm def}{=} -\sum_k \vartheta_k \pmod{2\pi}
\label{defect angle}
\ee
(subscript $k$ numbers the angles). Clearly, all deficit angles in our
construction are so far zero: the cluster of all 4-simplices having a given
common two-dimensional face can be mapped in the Euclidean space so that every
separate 4-simplex is mapped isometrically and the angles are taken with
the proper signs. We are going, however, to consider two kinds of
deformations of Euclidean values attributed to every 4-simplex that bring
about nonzero deficit angles.

The first kind of deformations consists in small changes of edge lengths:
we vary each of them arbitrarily in a neighbourhood of the value
obtained from the above construction. Every separate 4-simplex remains
Euclidean, but a mapping of their cluster in the Euclidean space
of the above type may no longer exist. As in Section~\ref{sec local},
we prefer to deal with {\em squared\/} edge lengths, denoting, for instance,
the squared length of edge~$AB$ as~$L_{AB}$. We use letters $a,b,\ldots$
for numbers of edges of the complex, and letters $i,j,\ldots$ --- for
numbers of two-dimensional faces. An example of these notations is the matrix of partial
derivatives $(\partial\omega_i/\,\partial L_a)$ which will play one of the main
r\^oles in our work.

The second kind of deformations of Euclidean quantities looks, probably,
more unusual. As is known, the {\em variations\/} of all metric values
in a given 4-simplex can be calculated from variations of the {\em areas\/}
of its ten {\em two-dimensional faces\/} (although a ``global'' statement
of such kind does not hold: areas, in contrast with lengths, do not
determine a 4-simplex uniquely). We take the values of areas obtained from
the above construction as initial ones and then give them arbitrary
independent variations. These variations determine the variations of
dihedral angles in each 4-simplex and, consequently, the variations of
deficit angles. In doing so, we {\em ignore the fact that the variations of
the length of one and the same edge may not coincide if we calculate
them from two neighbouring 4-simplices}.

Besides the angles at two-dimensional faces, we will need the angles at
{\em edges\/}. Consider at first one single 4-simplex. We define the angle
at its edge~$a$ as
\be
\Theta_a \stackrel{\rm def}{=} \sum_i \pa{S_i}{L_a} \,\vartheta_i,
\label{Theta}
\ee
where $i$ numbers the two-dimensional faces of the simplex (there are four
nonzero terms in the sum), while $S_i$ is, of course, the area of the
$i$th face. The reasonableness of such definition can be seen,
by the way, from the fact that thus defined angles obey an analogue
of the {\em Schl\"afli differential identity\/} (see~\cite{milnor}).

Recall that the usual Schl\"afli identity for a 4-simplex states that
\be
\sum_i S_i \,d\vartheta_i = 0
\label{schlafli}
\ee
for any infinitesimal deformations. We are going to prove that
$$
\sum_a L_a \,d\Theta_a = 0
$$
as well. To do so, it is enough to consider instead of $d\Theta_a$ the
partial derivatives w.r.t.\ an arbitrary~$L_b$:
\be
\sum_a L_a \pa{\Theta_a}{L_b} = \sum_a L_a \sum_i {\partial^2 S_i \over
\partial L_a \, \partial L_b}\, \vartheta_i + \sum_a L_a \sum_i \pa{S_i}{L_a}
\pa{\vartheta_i}{L_b}.
\label{proof modschlafli}
\ee
Now we note that $S_i$ are {\em homogeneous functions\/} of degree~$1$
in squared lengths~$L_a$. This yields $\sum_a L_a \,\partial S_i / \partial L_a
=S_i$ and $\sum_a L_a \,\partial^2 S_i / (\partial L_a \,\partial L_b) =0$.
Now we see that both terms in (\ref{proof modschlafli}) vanish
(the second one --- because of the usual Schl\"afli identity~(\ref{schlafli})).

To conclude this Section, we give a natural definition of the
{\em deficit angles around edges\/}:
\be
\Omega_a \stackrel{\rm def}{=} - \sum_k (\Theta_a)_k \pmod{2\pi},
\label{Omega}
\ee
where $k$ numbers 4-simplices that contain edge~$a$.

\section{Matrices relating various differentials: symmetry properties}
\label{sec matrices}

We now consider matrices that express the linear dependencies of the
differentials of deficit angles of two types introduced in
Section~\ref{sec differentials} on the differentials of squared lengths
or areas. These latter differentials, as was already explained, are regarded
as independent variables. Everything takes place in an infinitesimal
neighbourhood of the flat case (where all deficit angles are zero).

\begin{theorem}
The matrix $(\partial \omega_i / \partial S_j)$ is symmetric:
\be
\pa{\omega_i}{S_j} = \pa{\omega_j}{S_i}.
\label{sym}
\ee
\end{theorem}

{\it Proof}. Consider first a similar statement for one 4-simplex:
\be
\partial \vartheta_i / \partial S_j=\partial \vartheta_j /
\partial S_i.
\label{local sym}
\ee
To prove it, we consider a function
$f(S_1,\ldots,S_{10}) = \sum_i S_i \vartheta_i$, where the sum is taken over
all two-dimensional faces of the 4-simplex. The relation~(\ref{local sym})
follows from the equalness of the mixed derivatives
$\partial^2 f / (\partial S_i \, \partial S_j)=
\partial^2 f / (\partial S_j \, \partial S_i)$ if we use the
Schl\"afli formula~(\ref{schlafli}) while calculating them.
Now formula~(\ref{sym}) is obtained by summing over the dihedral angles
entering in the given deficit angle.

\medskip

Consider now the rectangular matrix $(\partial \Omega_a
/ \partial S_i)$. Note that the elements of this matrix are determined
correctly in spite of the fact that the definition of $\Omega_a$
involves some quantities (partial derivatives of the areas w.r.t.\ the
lengths) whose variations are not determined uniquely (they depend on
the given 4-simplex, like the variations of lengths do). The point is that
all such quantities are multiplied by the vanishing deficit angles~$\omega_j$ when we
calculate $\partial \Omega_a / \partial S_i$.

\begin{theorem}
\label{th conj}
Matrices $(\partial \Omega_a / \partial S_i)$ and
$(\partial \omega_i / \partial L_a)$ are mutually conjugate:
\be
\pa{\Omega_a}{S_i} = \pa{\omega_i}{L_a}.
\label{conjugate}
\ee
\end{theorem}

{\it Proof\/} follows from the chain of equalities:
\begin{eqnarray*}
\pa{\Omega_a}{S_i} = \pa{}{S_i} \left( \sum_j \pa{S_j}{L_a} \, \omega_j
\right) = \\
= \sum_j \pa{\omega_j}{S_i} \pa{S_j}{L_a} =
\sum_j \pa{\omega_i}{S_j} \pa{S_j}{L_a} = \pa{\omega_i}{L_a}.
\end{eqnarray*}
Here we have used, first, the expression of $\Omega_a$ in terms of $\omega_j$
(obtained from (\ref{defect angle}), (\ref{Theta}) and
(\ref{Omega})), secondly, the fact that the derivatives are taken with all
$\omega_j=0$ and, finally, the equality~(\ref{sym}).

\section{The invariant of moves $3\to 3$}
\label{sec invariant 3-3}

The preparatory work made in Sections \ref{sec differentials}
and~\ref{sec matrices} allows us now to pass from the local
formula (\ref{6term}) to its global analogue.

\begin{theorem}
\label{th restricted}
Choose in matrix $(\partial \omega_i/ \partial L_a)$ a largest
square submatrix~$\cal B$ with nonzero determinant. Let
$\cal B$ contain a row corresponding to such face $i=ABC$ that
belongs to exactly three 4-simplices. Then those latter can be replaced
by three new 4-simplices, as in Section~\ref{sec local}. After such
replacement, take in the\/
{\em new\/} matrix $(\partial \omega_i/ \partial L_a)$ the new
submatrix~$\cal B$ containing the same rows and columns, with the only
following change: the row corresponding to face~$ABC$ is replaced by the
row corresponding to the new face~$DEF$.

The expression
\be
{\displaystyle \det {\cal B}\cdot \prod_{
\hbox{\scriptsize
\begin{tabular}{c}
\rm over all\\[-.5\baselineskip]
\rm 4-simplices
\end{tabular} } } V
\over \displaystyle \prod_{
\hbox{\scriptsize
\begin{tabular}{c}
\rm over all\\[-.5\baselineskip]
\rm 2-dim.\ faces
\end{tabular} } } S }
\label{3-3 restricted}
\ee
does not change under such re-building.
\end{theorem}

{\it Proof}. The only thing in formula (\ref{3-3 restricted}) which is not
completely evident when compared with (\ref{6term}) is that there is
a determinant $\det {\cal B}$ instead of just one
differential~$d\omega$. It is convenient to write that determinant as
$\displaystyle \det {\cal B} = \bigwedge d\omega \,/\,\bigwedge dL \,$,
where the exterior product in the numerator is taken over the two-dimensional
faces corresponding to the rows of matrix~$\cal B$, while the one in the
denominator is taken over the edges corresponding to the columns of~$\cal B$.
Now we note that if $\omega_{ABC}=0$ or, equivalently, $\omega_{DEF}=0$, then
the deficit angles at the {\em remaining\/} faces do not change under
the re-building (because they are obtained from the same dihedral angles).
Thus, the difference $d\omega_j - d\omega'_j$ for any of the remaining faces
(with those differentials considered as depending on~$dL$'s),
taken before and after the re-building, is proportional to
$d\omega_{ABC}$ or, equivalently, to~$d\omega_{DEF}$.
Consequently, $\det {\cal B}$ gets multiplied by
$d\omega_{ABC} / \,d\omega_{DEF}$ under the re-building. The theorem is proven.

\medskip

Expression (\ref{3-3 restricted}) depends on a specific choice of the
subsets of faces and edges corresponding to the rows and columns of
matrix~$\cal B$. In getting rid of these dependencies, we will be helped,
like in the three-dimensional case of papers~\cite{3dim1,3dim2},
by some {\em differential forms}.

We denote by~$\cal D$ the set of two-dimensional faces corresponding to the
rows of matrix~$\cal B$. By $\ol{\cal D}$ we denote the set of the
{\em rest\/} of the two-dimensional faces in the complex. Similarly,
the set of edges corresponding to the columns of~$\cal B$ will be
denoted~$\cal C$, while $\ol{\cal C}$ will be the set of the remaining
edges in the complex. We write, somewhat loosely:
$$
{\cal B} =
\left.\vphantom{\pa{}{}}\right._{\cal D}\!
\left( \pa{\omega_i}{L_a} \right)
\!\!\left.\vphantom{\pa{}{}}\right._{\cal C}.
$$
The looseness consists in the fact that, for matrix~$\cal B$ to be
determined completely, we need also the {\em ordering\/} of rows and columns,
and this ordering influences the sign of determinant~$\det {\cal B}$.

At the moment, we assume that the set ${\cal D} \cup \ol{\cal D}$ of all
two-dimensional faces is fixed. The set ${\cal C} \cup \ol{\cal C}$ is fixed,
too, because the set of all edges does not change under re-buildings
$3\to 3$ which we are considering.

\begin{theorem}
\label{th invCD}
The expression
\be
(\det {\cal B})^{-1} \cdot \bigwedge_{i\in \ol{\cal D}} dS_i \cdot
\bigwedge_{a\in \ol{\cal C}} dL_a\, ,
\label{invCD}
\ee
taken to within its sign, does not depend on the choice of sets $\cal C$
and~$\cal D$.
\end{theorem}

We explain the meaning of expression~(\ref{invCD}). Set~$\ol{\cal C}$ is,
in full analogy with papers~\cite{3dim1,3dim2}, a largest subset in the
set of all edges such that the lengths of edges from~$\ol{\cal C}$ can be
given arbitrary infinitesimal increments $dL_a$ {\em without violating
the ``zero curvature'' conditions\/} $d\omega_i=0$ for all~$i$. The rest of
$dL_a$ (for $a\in {\cal C}$) are determined uniquely from those conditions.
This allows us to express
$\displaystyle \bigwedge_{a\in \ol{\cal C}_1} dL_a $
in terms of
$\displaystyle \bigwedge_{a\in \ol{\cal C}} dL_a $,
where $\ol{\cal C}_1$ is another set similar to~$\ol{\cal C}$. To be exact,
these exterior products are proportional, and we can calculate the
proportionality factor from matrix~$(\partial \omega_i /\, \partial L_a)$.
This means that we consider exterior products of such kind
as products of linear forms defined on the linear
space of such infinitesimal deformations of lengths which leave the
deficit angles zero.

It turns out that set~$\ol {\cal D}$ admits a similar description:
$\ol {\cal D}$ is a largest subset in the set of all two-dimensional faces
in the complex such that the areas of faces in $\ol {\cal D}$ can be given
arbitrary infinitesimal increments $dS_i$ without violating the
{\em ``zero curvature around edges''\/} conditions $d\Omega_a = 0$ for
all~$a$. Here we regard $d\Omega_a$ as functions of all~$dS_i$, as in
Section~\ref{sec matrices}. The values $dS_i$ for $i\in {\cal D}$
are determined from the conditions $d\Omega_a = 0$ unambiguously.

Indeed, Theorem~\ref{th conj} states that matrices
$(\partial \Omega_a / \, \partial S_i)$ and
$(\partial \omega_i / \, \partial L_a)$ are obtained one from the other
by matrix transposing (which we denote by superscript~$\rm T$).
Consequently, ${\cal B}^{\rm T} =
\left.\vphantom{|}\right._{\cal C}\!
(\partial \Omega_a / \, \partial S_i )
\!\!\left.\vphantom{|}\right._{\cal D}$
is a maximal nondegenerate square submatrix in
$(\partial \Omega_a / \, \partial S_i)$, and this yields our statements.

Thus, we consider exterior products of area differentials
on the linear space of such infinitesimal deformations of
the areas that leave the deficit angles {\em around edges\/} zero. If
$\ol{\cal D}_1$ is a different set of two-dimensional faces similar to
$\ol{\cal D}$ then the forms
$\displaystyle \bigwedge_{i\in \ol {\cal D}_1} dS_i$ and
$\displaystyle \bigwedge_{i\in \ol {\cal D}} dS_i$ are proportional, with
the proportionality factor calculated again from
matrix~$(\partial \omega_i /\, \partial L_a)$.

\medskip

{\it Proof of Theorem~\ref{th invCD}}. We prove first that expression
(\ref{invCD}) does not change when $\cal C$ and~$\ol{\cal C}$ exchange their
elements. For simplicity, we limit ourselves to the case where
only {\em one\/} edge~$b$ is carried from $\ol{\cal C}$ to~$\cal C$,
while an edge~$c$ takes its place. In doing so,
$\displaystyle \bigwedge_{a\in \ol {\cal C}} dL_a$ is of course multiplied by
\be
\left.\pa{L_c}{L_b}\right|_{\hbox{\small \rm with all $d\omega_i = 0$
and $dL_a=0$, where $a\in \ol{\cal C},\; a\ne b$}}\; .
\label{dL/dL}
\ee

On the other hand, we can calculate the partial derivative (\ref{dL/dL})
by considering the squared edge lengths in~$\cal C$ as {\em implicit functions\/}
of the squared edge lengths in~$\ol {\cal C}$ given by the conditions
$\omega_i=0$, $i\in {\cal C}$. This yields
$$
\pa{L_c}{L_b} = - {\det {\cal B}_{\rm new} \over \det {\cal B}}\,.
$$
Consequently, $\det {\cal B}$ gets multiplied, up to a sign, by the same
expression~(\ref{dL/dL}) under such exchange~$b\leftrightarrow c$.

We prove now that expression (\ref{invCD}) does not change when
$\cal D$ and~$\ol{\cal D}$ exchange their elements. Again, we assume for
simplicity that only one pair of elements is exchanged. To be exact, let there be
$n$~two-dimensional faces in the complex, and let set~$\cal D$ include
the faces with numbers from $1$ through $m$, while $\ol{\cal D}$ --- 
those with numbers from $m+1$ through $n$. We are going to consider the
interchange~$m\leftrightarrow m+1$.

Certainly, $m=\rank (\partial\omega_i / \, \partial L_a)$. Thus, in
particular, $(m+1)$th row of matrix $(\partial\omega_i / \, \partial L_a)$
is a linear combination of the first $m$ rows with some coefficients
$c_1,\ldots,\allowbreak c_m$. This can be written as the relation
\be
d\omega_{m+1} = c_1 \,d\omega_1 + \cdots + c_m \,d\omega_m,
\label{domega}
\ee
and also as
\be
\pmatrix{ -c_1 & \ldots & -c_m & 1 & 0 & \ldots & 0 } \left(
\pa{\omega_i}{L_a}\right) =0.
\label{c1cm}
\ee

Matrix transposing relation (\ref{c1cm}), we find
$$
\left( \pa{\Omega_a}{S_i} \right)
\pmatrix{ -c_1 \cr \ldots \cr -c_m \cr 1 \cr 0 \cr \vdots \cr 0 } =0.
$$
This means that all $d\Omega_a=0$ as soon as the following proportions
are satisfied:
\begin{eqnarray*}
dS_1:dS_2:\cdots:dS_m:dS_{m+1}:dS_{m+2}:\cdots:dS_n= \\
= -c_1:-c_2:\cdots:-c_m:1:0:\cdots:0.
\end{eqnarray*}

In particular,
\be
\left.\pa{S_m}{S_{m+1}}\right|_{\hbox{\small\rm with all $d\Omega_a=0$
and $dS_{m+2}=\cdots=dS_n=0$}} = -c_m.
\label{dSdS}
\ee
Note that $\displaystyle \bigwedge_{i\in \ol{\cal D}} dS_i$
gets multiplied exactly by (\ref{dSdS}) under the interchange of faces
$m \leftrightarrow m+1$.

On the other hand, equality (\ref{domega}) shows that $\det {\cal B}$
gets multiplied by~$c_m$ under such interchange of faces. Thus,
Theorem~\ref{th invCD} is proven.

\medskip

Combining Theorem~\ref{th invCD} with Theorem~\ref{th restricted}, we get the
desired invariant of moves~$3\to 3$.

\begin{theorem}
\label{th inv}
The value
\be
\displaystyle \prod_{
\hbox{\scriptsize
\begin{tabular}{c}
\rm over all\\[-.5\baselineskip]
\rm 2-dim.\ faces
\end{tabular} } } S
\cdot
\bigwedge_{\ol {\cal D}} dS
\cdot
\bigwedge_{\ol {\cal C}} dL
\over
\det {\cal B}
\cdot
\displaystyle \prod_{
\hbox{\scriptsize
\begin{tabular}{c}
\rm over all\\[-.5\baselineskip]
\rm 4-simplices
\end{tabular} } } V
\label{inv}
\ee
does not depend on the choice of sets $\cal C$ and $\cal D$ and is an
invariant of moves~$3\to 3$.
\end{theorem}

\section{Discussion}
\label{sec discussion}

It is not hard to see in the expression~(\ref{inv}) a direct generalization
of the similar invariant of three-dimensional manifolds.
To be exact, (\ref{inv}) looks very much like the squared expression
for the three-dimensional manifold invariant from paper~\cite{3dim1}
(see~\cite{3dim1}, formula~(30)). Note that the mentioned three-dimensional
manifold invariant is an invariant of {\em all\/} kinds of re-buildings of
three-dimensional simplicial complexes, or Pachner moves, namely,
$2\leftrightarrow 3$ and $1\leftrightarrow 4$. This makes us think that in
the four-dimensional case, too, the applicability of expression~(\ref{inv})
goes beyond the limits of just moves~$3\to 3$. Besides, this makes us hope
for nontrivial consequences of Theorem~\ref{th inv}, because the mentioned
three-dimensional invariant passes the nontriviality test~\cite{3dim2}.

In the further works we plan to study the behaviour of expression~(\ref{inv})
under re-buildings $2\leftrightarrow 4$ and~$1\leftrightarrow 5$.
Meanwhile, we propose the following argument in defence of investigating
the moves $3\to 3$ even by themselves, by referring to the analogy with the
{\em Yang--Baxter equation\/} (YBE). As is
well-known, YBE corresponds, from the knot theory viewpoint,
to only one of the {\em Reidemeister moves\/} which are fundamental there.
Nevertheless, YBE proves to be interesting enough by itself from an algebraic
point of view, and finds nontrivial applications, in particular, in
statistical physics and quantum field theory.

\medskip

{\bf Acknowledgements. }\vadjust{\nobreak}I am grateful to John Barrett and
Justin Roberts for their help on the issues related to the Schl\"afli
differential identity. The work has been performed with a partial financial
support from Russian Foundation for Basic Research under
Grant no.~01-01-00059.

\end{document}